\documentclass[a4paper]{article}
\usepackage{RR}
\usepackage{hyperref}

\usepackage{epic,eepic,epsfig,amssymb,amsmath,amsthm}
\newtheorem{theorem}{Theorem}

\newtheorem{lemma}{Lemma}[section]
\theoremstyle{remark}\newtheorem{remark}{Remark}[section]
\theoremstyle{definition}\newtheorem{definition}{Definition}[section]

\theoremstyle{definition}

\def\R{\textrm{I\kern-0.21emR}}
\def\N{\textrm{I\kern-0.21emN}}

\RRdate{Mars 2007}

\RRauthor{

Fr\'ed\'eric Bonnans\thanks[sfn]{CMAP UMR CNRS 7641, Ecole Polytechnique, 
and INRIA Futurs, 91128 Palaiseau, France} 
\thanks{\texttt{Frederic.Bonnans@polytechnique.edu}}

\and

Pierre Martinon\thanksref{sfn}
\thanks{\texttt{martinon@cmap.polytechnique.fr}}

\and

Emmanuel Tr\'elat\thanks{Universit\'e d'Orl\'eans, Math., Labo. MAPMO, UMR CNRS 6628, Route de Chartres, BP 6759, 45067 Orl\'eans cedex 2, France (\texttt{Emmanuel.Trelat@univ-orleans.fr})}
\thanksref{sfn}

}
\authorhead{Bonnans, Martinon \& Tr\'elat}
\RRtitle{Arcs singuliers dans le probl\`eme de Goddard g\'en\'eralis\'e}
\RRetitle{Singular arcs in the generalized Goddard's Problem}

\titlehead{Example of RR.sty}
\RRnote{Support from the French space agency CNES (Centre National d'Etudes Spatial) is gratefully acknowledged.}

\RRresume{On e\'tudie des variantes du probl\`eme de Goddard pour des trajectoires non verticales. Le contr\^ole est la force de pouss\'ee, et l'objectif est de maximiser un certain crit\`ere, typiquement la masse finale.
Dans ce rapport, on prouve par une analyse bas\'ee sur le Principe du Maximum de Pontryagin que les trajectoires optimales peuvent comporter des arcs singuliers (sur lesquels le module de la pouss\'ee n'est ni nul ni maximal), que l'on peut caract\'eriser et calculer. On pr\'esente des simulations num\'eriques, \`a la fois avec des m\'ethodes directes et indirectes, qui montrent la pertinence de la prise en compte des arcs singuliers dans la strat\'egie de contr\^ole. La m\'ethode indirecte utilis\'ee est bas\'ee sur l'analyse th\'eorique mentionn\'ee pr\'ec\'edemment, et consiste \`a combiner une m\'ethode de tir avec une approche homotopique. L'homotopie prend la forme d'une r\'egularisation quadratique du probl\`eme, et permet de traiter l'aspect non lisse du contr\^ole optimal.}

\RRabstract{We investigate variants of Goddard's problems for nonvertical
trajectories. The control is the thrust force,
and the objective is to maximize a certain final cost, typically, the final mass.
In this report, performing an analysis based on the Pontryagin Maximum Principle, we prove that optimal trajectories may involve singular arcs (along which the norm of the thrust is neither zero nor maximal), that are computed and characterized. Numerical simulations are carried out, both with direct and indirect methods, demonstrating the relevance of taking into account singular arcs in the control strategy. The indirect method we use is based on our previous theoretical analysis and consists in combining a shooting method with an homotopic method. The homotopic approach leads to a quadratic regularization of the problem and is a way to tackle with the problem of nonsmoothness of the optimal control.}

\RRmotcle{Contr\^ole optimal, probl\`eme de Goddard, trajectoires singuli\`eres, m\'ethode de tir, homotopie, m\'ethodes directes.}
\RRkeyword{Optimal control, Goddard's problem, 
singular trajectories, shooting method, homotopy, direct methods.}

\RRprojets{Commands}
\RRtheme{\THNum} 

\URFuturs 

\begin{document}
\makeRR   

\section{Introduction}
The classical Goddard's problem (see \cite{Goddard,SeyCli93,TsiKel91})
consists in maximizing the final altitude of a rocket with vertical
trajectory, the controls being the norm and direction of the thrust force.
Due to nonlinear effects of aerodynamic forces, the optimal strategy may involve subarcs along which the thrust is neither zero nor equal to its maximal value, namely, since the control variable enters linearly in the dynamics and the cost function is over the final cost, \textit{singular arcs}. A natural extension of this model for nonvertical trajectories is the control system
\begin{equation}\label{system}
\begin{split}
\dot{r}&=v ,\\
\dot{v}&=-\frac{D(r,v)}{m}\frac{v}{\Vert v\Vert} - g(r) + C\frac{u}{m},\\
\dot{m}&=-b\Vert u\Vert,
\end{split}
\end{equation}
where the state variables are $r(t)\in\R^3$ (position of the spacecraft), $v(t)\in\R^3$ (velocity vector) and $m(t)$ (mass of the engine).
Also, $D(r,v)>0$ is the drag component, $g(r)\in\R^3$ is the usual gravity force, and $b$ is a positive real number depending on the engine. The thrust force is $Cu(t)$, where
$C>0$ is the maximal thrust, and the control is the normalized thrust
$u(t)\in\R^3$, submitted to the constraint
\begin{equation}\label{constraint}
\Vert u(t)\Vert\leq 1.
\end{equation}
The real number $b>0$ is such that the speed of ejection is $C/b$. 
Here, and throughout the paper, $\Vert\ \Vert$ denotes the usual Euclidean norm in $\R^3$.

We consider the optimal control problem of steering the system from a given initial point
\begin{equation}\label{initialconditions}
r(0)=r_0,\ v(0)=v_0,\ m(0)=m_0,
\end{equation}
to a certain target $M_1\subset\R^7$, in time $t_f$ that may be fixed or not, while maximizing a final cost. For the moment, there is no need to be more specific with final conditions and the cost.
In real applications, the problem is typically to reach a given orbit, either in
minimal time with a constraint on the final mass, or by maximizing
the final mass, or a compromise between the final mass and time
to reach the orbit.
In our numerical experiments we will study the problem of
maximizing the final mass (i.e., minimizing the fuel consumption) subject to a fixed final position $r(t_f)=r_f$, the final velocity vector and final time being free.

Depending on the features of the problem (initial and final conditions, mass/thrust ratio, etc), it is known that control strategies that consist in choosing the control so that $\Vert u(t)\Vert$ is piecewise constant all along the flight, either equal to $0$ or to the maximal authorized value $1$, may not be optimal, as a consequence of the high values of the drag for high speed. Optimal trajectories may indeed involve singular arcs, and it is precisely the aim to this article to perform such an analysis and prove that the use of singular arcs is relevant in the problem of launchers.

\medskip

The article is structured as follows. In Section \ref{secprelim}, we recall the Pontryagin Maximum Principle, and the concept of singular trajectories. A precise analysis of the optimal control problem is performed in Section \ref{secanalysis}, where extremals are derived, and singular trajectories are computed. Theorem \ref{thmcontroles} makes precise the structure of the optimal trajectories. Section \ref{secnum} is devoted to numerical simulations. The problem is first implemented with indirect methods, based on our theoretical analysis with the maximum principle, and, numerically, our method uses a shooting method combined with an homotopic approach. The homotopic method, leading to a quadratic regularization, permits to tackle with the problem of nonsmoothness of the optimal control. Experiments are also made using direct methods, i.e., by discretizing control variables and solving the resulting nonlinear optimization problem. Less precise than the indirect one, this method permits however to validate our approach by checking that results are consistent with the previously computed solution.

Our results show, as expected, that taking into account singular arcs in the control strategy permits to improve slightly the optimization criterion. The numerical simulations presented in this paper, using a simplified and more academic model and set of parameters, constitute the first step in the study of a realistic launcher problem. 


\section{Preliminaries}\label{secprelim}
In this section we recall a general version of the Pontryagin Maximum Principle (see \cite{PBGM}, and for instance \cite{BryHo75} for its practical application), and a definition and characterizations of singular arcs.

Consider the autonomous control system in $\R^n$
\begin{equation}\label{systemPMP}
\dot{x}(t)=f(x(t),u(t)),
\end{equation}
where $f:\R\times\R^n\times\R^m\longrightarrow\R^n$ is of class $C^1$, and where the controls are measurable and bounded functions defined on a subinterval $[0,t_e(u)[$ of $\R^+$ with values in $\Omega \subset \R^m$. Let $M_0$ and $M_1$ be subsets of $\R^n$. Denote by $\mathcal{U}$ the set of admissible controls $u$, whose associated trajectories are well defined and join an initial point in $M_0$ to a final point in $M_1$, in time $t(u)<t_e(u)$.

Define the cost of a control $u$ on $[0,t]$ by
$$C(t,u)=\int_0^{t}f^0(x(s),u(s))ds +g^0(t,x(t)), $$
where $f^0:\R^n\times\R^m\longrightarrow\R$ and $g^0:\R^n\rightarrow\R$ are of class $C^1$, and $x(\cdot)$ is the trajectory solution of (\ref{system}) associated to the control $u$.

Consider the optimal control problem of finding a trajectory joining
$M_0$ to $M_1$ and minimizing the cost. The final time may be free or not.


\subsection{Pontryagin Maximum Principle}\label{secPMP}
According to the Pontryagin Maximum Principle (see \cite{PBGM}),
if the control $u\in{\mathcal{U}}$ associated to the trajectory $x(\cdot)$
is optimal on $[0,T]$, then there exists an absolutely continuous mapping
$p(\cdot): [0,T]\longrightarrow \R^n$ called \textit{adjoint vector}, and a real number
$p^0\leq 0$, such that the couple $(p(\cdot),p^0)$ is nontrivial, and such that, for almost every $t\in[0,T]$,
\begin{equation} \label{systPMP}
\begin{split}
\dot{x}(t)&=\frac{\partial H}{\partial p}(x(t),p(t),p^0,u(t)),\\
\dot{p}(t)&=-\frac{\partial H}{\partial x}(x(t),p(t),p^0,u(t)),
\end{split}
\end{equation}
where $H(x,p,p^0,u)=\langle p,f(x,u)\rangle+p^0f^0(x,u)$ is the \textit{Hamiltonian} of the optimal control problem. Moreover, the function
$$t \longmapsto \max_{v\in\Omega}
H(x(t),p(t),p^0,v)$$
is constant on $[0,T]$,
and the \textit{maximization condition}
\begin{equation} \label{contraintePMP}
H(x(t),p(t),p^0,u(t))=\max_{v\in\Omega} H(x(t),p(t),p^0,v)
\end{equation}
holds almost everywhere on $[0,T]$.

Moreover, if the final time $T$ to join the target set $M_1$ is free, then
\begin{equation}\label{Hamnul}
\max_{v\in\Omega}
H(x(t),p(t),p^0,v) = -p^0\frac{\partial g^0}{\partial t}(T,x(T)).
\end{equation}
for every $t\in[0,T]$.

Furthermore, if $M_0$ and $M_1$ (or just one of them) are submanifolds of $\R^n$ having tangent spaces in $x(0)\in M_0$ and
$x(T)\in M_1$, then the adjoint vector can be chosen so as to satisfy the \textit{transversality conditions} at both extremities (or just one of them)
\begin{equation}\label{condt1}
p(0)\ \bot\ T_{x(0)}M_0
\end{equation}
and
\begin{equation}\label{condt2}
p(T)-p^0\frac{\partial g^0}{\partial x}(T,x(T))\ \bot\ T_{x(T)}M_1.
\end{equation}

\medskip

\noindent
An \textit{extremal} of the optimal control problem is a fourth-tuple
$(x(\cdot),p(\cdot),p^0,u(\cdot))$ solution of (\ref{systPMP}) and
(\ref{contraintePMP}). If $p_0=0$, then the extremal is said to be
{\textit{abnormal}}, and if $p^0\neq 0$ then the extremal is said to be
{\textit{normal}}.


\subsection{Singular arcs} \label{secchar}
Given $x_0\in\R^n$ and two real numbers $t_0$, $t_1$, with $t_0<t_1$, denote by ${\cal U}_{x_0,t_0,t_1}$ the set of controls $u\in L^\infty([t_0,t_1],\Omega1)$, with $\Omega1$ an open subset of $\Omega$, such that the trajectory
$t\mapsto x(t,x_0,t_0,u)$, solution of (\ref{system}), associated with the control $u$ on $[t_0,t_1]$, and such that $x(t_0)=x_0$, is well defined on $[t_0,t_1]$.
Define the \textit{end-point mapping} $E_{x_0,t_0,t_1}$ by
$E_{x_0,t_0,t_1}(u) := x(t_1,x_0,t_0,u),$
for every $u\in{\cal U}_{x_0,t_0,t_1}$.
It is classical that $E_{x_0,t_0,t_1}:{\cal U}_{x_0,t_0,t_1} \rightarrow
\R^n$ is a smooth map.

A control $u\in {\cal U}_{x_0,t_0,t_1}$ is said to be \textit{singular} if $u$
is a critical point of the end-point mapping $E_{x_0,t_0,t_1}$,
i.e., its differential $dE_{x_0,t_0,t_1}(u)$ at $u$ is not surjective.
In this case, the trajectory $x(\cdot,x_0,t_0,u)$ is said to be \textit{singular} on $[t_0,t_1]$.

\medskip

Recall the two following standard characterizations of singular controls (see \cite{BonChy03,PBGM}). A control $u\in {\cal U}_{x_0,t_0,t_1}$ is singular if and only if the
linearized system along the trajectory $x(\cdot,x_0,t_0,u)$ on $[t_0,t_1]$ is not
controllable. This is also equivalent to the existence of an absolutely continuous mapping
$p_1:[t_0,t_1]\longrightarrow \R^n \setminus \{0\}$ such that, for almost every $t\in[t_0,t_1]$,
$$\dot{x}(t)=\frac{\partial H_1}{\partial p}(x(t),p_1(t),u(t)) ,  \quad
\dot{p}(t)=-\frac{\partial H_1}{\partial x}(t,x(t),p_1(t),u(t)), $$
$$\frac{\partial H_1}{\partial u}(x(t),p_1(t),u(t))=0 ,$$
where $H_1(x,p_1,u)=\langle p_1,f(x,u)\rangle$ is the Hamiltonian of the system.

Note that singular trajectories coincide with projections of abnormal extremals for which the maximization condition (\ref{contraintePMP}) reduces to $\frac{\partial H}{\partial u}=0$.

\medskip

For a given trajectory $x(\cdot)$ of the system (\ref{system}) on $[0,T]$, associated to a control $u\in {\cal U}_{x(0),0,T}$, we say that $x(\cdot)$ involves a singular arc, defined on the subinterval $[t_0,t_1]\subset[0,T]$, whenever the control $u_{\vert [t_0,t_1]}$ for the control system restricted to $[t_0,t_1]$ is singular.

In the case when the dynamics $f$ and the instantaneous cost $f^0$ are linear in the control $u$, a singular arc corresponds to an arc along which one is unable to compute the control \textit{directly} from the maximization condition of the Pontryagin maximum principle (at the contrary of the bang-bang situation). Indeed, in this case, the above condition $\frac{\partial H_1}{\partial u}=0$ along the arc means that some function (called switching function) vanishes identically along the arc. Then, it is well known that, in order to derive an expression of the control along such an arc, one has to differentiate this relation until the control appears explicitly.
It is also well known that such singular arcs, whenever they occur, may be optimal. Their optimal status may be proved using generalized Legendre-Clebsch type conditions or the theory of conjugate points (see \cite{Rob67,Goh}, or see \cite{AgSachkov,BCT} for a complete second-order optimality theory of singular arcs).


\section{Analysis of the optimal control problem}\label{secanalysis}
With respect to the notations used in the previous section, we set 
$$x=\begin{pmatrix}
r\\ v\\ m
\end{pmatrix}\in \R^3\times\R^3\times\R,\quad
f(x,u)=\begin{pmatrix}
v \\
-\frac{D(r,v)}{m}\frac{v}{\Vert v\Vert} - g(r) + C\frac{u}{m}\\
-b\Vert u\Vert
\end{pmatrix},
$$
and $f^0=0$. Here, the set $\Omega$ of constraints on the control is the closed unit ball of $\R^3$, centered at $0$.

Consider the optimal control problem of minimizing some final cost $g^0(t_f,x(t_f))$, for the control system (\ref{system}), with initial conditions (\ref{initialconditions}) and final conditions $x(t_f)\in M_1$ in time $t_f$ which may be free or not.

We make the following assumption.

\medskip

\noindent{\bf Assumption $\bf (H)$.} The function $g^0$ is such that:
\begin{itemize}
\item either the final mass $m(t_f)$ is free, and $\frac{\partial g^0}{\partial m}\neq 0$,
\item or the final time $t_f$ is free, and $\frac{\partial g^0}{\partial t}\neq 0$.
\end{itemize}

In the first situation, the target set $M_1\subset\R^7$ can be written as $M_1=N_1\times\R$, where $N_1$ is a subset of $\R^6$.
A typical example is the problem of maximizing the final mass, for which $g^0(t,x)=-m$. If the final condition is $r(t_f)=r_1$ and $\Vert v(t_f)\Vert=a$, then $M_1=\{r_1\}\times S(0,a)\times\R$, where $S(0,a)$ is the sphere of $\R^3$, centered at $0$, with radius $a$.

In the second situation, a typical example is the minimal time problem to reach some target. In this case, $g^0(t,x)=t$.


\subsection{Computation of extremals}
According to Section \ref{secPMP}, the Hamiltonian of the optimal control problem under consideration is
\begin{equation}\label{defH}
H=\langle p_r,v\rangle + \left\langle p_v, -\frac{D(r,v)}{m}\frac{v}{\Vert v\Vert} - g(r) + C\frac{u}{m}  \right\rangle - p_m b \Vert u\Vert,
\end{equation}
where $\langle\ ,\ \rangle$ denotes the usual scalar product in $\R^3$. Here, the adjoint vector is denoted by
$$p(t)=\begin{pmatrix}
p_r(t)\\ p_v(t)\\ p_m(t)\end{pmatrix}\in \R^3\times\R^3\times\R .$$
In what follows, we assume the mappings $D$ and $g$ to be of class $C^1$.
Applying Pontryagin's Maximum Principle leads to the adjoint equations
\begin{equation}\label{systemeadjoint}
\begin{split}
\dot{p}_r &= \frac{1}{m}\frac{\langle p_v, v\rangle}{\Vert v\Vert} \frac{\partial D}{\partial r} + \left\langle p_v, \frac{\partial g}{\partial r} \right\rangle , \\
\dot{p}_v &= -p_r + \frac{1}{m}\frac{\langle p_v, v\rangle}{\Vert v\Vert} \frac{\partial D}{\partial v} + \frac{D}{m} \frac{ p_v}{\Vert v\Vert}  -\frac{D}{m}\langle p_v,v\rangle \frac{v}{\Vert v\Vert^3} , \\
\dot{p}_m &= \frac{1}{m} \left\langle p_v, -\frac{D(r,v)}{m}\frac{v}{\Vert v\Vert} + C\frac{u}{m} \right\rangle.
\end{split}
\end{equation}
Moreover, if $u$ is an optimal control on $[0,t_f]$, then, for almost every $t\in[0,t_f]$,
$u(t)$ maximizes the function
$$
\Phi_t(w) := \frac{C}{m(t)}\langle p_v(t),w\rangle - bp_m(t)\Vert w\Vert,
$$
among all possible $w\in\R^3$ such that $\Vert w\Vert \leq 1$.

The next technical lemma is the first step in the analysis of extremals.

\begin{lemma}\label{lempv}
If there exists $t_0\in [0,t_f]$ such that $p_r(t_0) = p_v(t_0) = 0$, then 
$p_r(t) = p_v(t) = 0$, and $p_m(t)=p_m(t_f)$, 
for every $t\in [0,t_f]$. Moreover, $p_m(t_f)\neq 0$, and if $p_m(t_f) >0$ then $u(t) = 0$
on $[0,t_f]$, otherwise $\|u(t)\| = 1$ on $[0,t_f]$. 
\end{lemma}

\begin{proof}
The first statement follows immediately from a uniqueness argument applied to the system
\eqref{systemeadjoint}. It follows from the expression of the Hamiltonian
function that, if $p_m(t) >0$, then $u(t) = 0$, and if $p_m(t) <0$, 
then $\|u(t)\| = 1$.
In the first case of Assumption $(H)$, the transversality condition (\ref{condt2}) yields in particular
$$
p_m(t_f) = p^0\frac{\partial g^0}{\partial m}(t_f,x(t_f)).
$$
Therefore, $p_m(t)$ cannot be equal to zero (otherwise the adjoint vector $(p,p^0)$
would be zero, contradicting the maximum principle).
In the second case of Assumption $(H)$, it follows from (\ref{Hamnul}) and (\ref{defH}) that
$$
p_m(t)b\Vert u(t)\Vert = p^0\frac{\partial g^0}{\partial t}(t_f,x(t_f)).
$$
Therefore, similarly, $p_m(t)$ cannot be equal to zero.
The conclusion follows.
\end{proof}

An extremal satisfying the conditions of Lemma \ref{lempv} (ie $p_r(t)=p_v(t)=0$ for every $t\in [0,t_f]$) is called \textit{degenerate}. For such extremals, the control is either identically equal to zero, or or maximal norm, along the whole trajectory. Such kind of trajectories can be excluded for practical applications and are thus discarded in the sequel.

\begin{lemma}\label{lembang}
Consider a nondegenerate extremal. Then:
\begin{enumerate}
\item The set $\mathcal{T} := \{ t\in[0,t_f]\ \vert\ p_v(t)= 0 \}$ has a finite cardinal.
\item There exists a measurable function $\alpha$ on $[0,t_f]$, with values in $[0,1]$, such that
\begin{equation}\label{ucolinpv}
u(t) = \alpha(t) \frac{ p_v(t) }{ \|p_v(t)\| } ,\quad\textrm{a.e. on}\ [0,t_f].
\end{equation}
\item
Set $\Psi(t) := \frac{C}{m(t)}\Vert  p_v(t)\Vert - b p_m(t)$. Then,
\begin{equation*}
\alpha(t) = \left\{ \begin{array}{ll}
0 & \textrm{if}\ \Psi(t)<0,\\
1 & \textrm{if}\ \Psi(t)>0.
\end{array}\right.
\end{equation*}
\end{enumerate}
\end{lemma}

\begin{proof}
If $t\in \cal T$, then by the costate equation \eqref{systemeadjoint},
$\dot p_v(t) = - p_r(t)$ is not zero (since the extremal is not degenerate). Therefore $\cal T$ has only isolated points, and hence, has a finite cardinal.

Writing $w=\alpha d$, with $\alpha=\|w\|$ and $d$ of unit norm, we get
$
\Phi_t(w) = \alpha \left(
\frac{C}{m(t)}\langle p_v(t),d\rangle - bp_m(t)\right).
$
Since $p_v(t) \neq 0$ a.e., points 2 and 3 of the lemma follow immediately from the maximization condition.
\end{proof}

The continuous function $\Psi$ defined in Lemma \ref{lembang} is called \textit{switching function}.
In the conditions of the lemma, the extremal control is either equal to $0$, or saturating the constraint and of direction $p_v(t)$. The remaining case, not treated in this lemma and analyzed next, is the case where the function $\Psi$ vanishes on a (closed) subset $I\subset [0,t_f]$ of positive measure.

\begin{remark}
Let $[t_0,t_1]$ be a subinterval of $I$ on which $\alpha(t)>0$. Then, the control $u_{\vert [t_0,t_1]}$ is singular.

Indeed, it suffices to notice that, using \eqref{ucolinpv},
$$\frac{\partial \Phi_t}{\partial w}(u(t))= \left(\frac{C}{m(t)}\Vert  p_v(t)\Vert - b p_m(t)\right)\frac{p_v(t)}{\Vert p_v(t)\Vert}=\Psi(t)\frac{p_v(t)}{\Vert p_v(t)\Vert} , $$
and to use the Hamiltonian characterization of singular controls recalled in Section \ref{secchar}.
\end{remark}

Singular arcs may thus occur in our problem whenever $\Psi$ vanishes, and we next provide an analysis of that case, and show how to derive an expression of such singular controls.


\subsection{Analysis of singular arcs}
Throughout this section, we assume that
\begin{equation}\label{rel1}
\Psi(t) = \frac{C}{m(t)}\Vert  p_v(t)\Vert - b p_m(t) = 0
\end{equation}
for every $t\in I$, where
$I$ is a (closed) measurable subset of $[0,t_f]$ of positive Lebesgue measure.

Usually, singular controls are computed by derivating this relation with respect to $t$, until $u$ appears explicitly. The following result is required (see \cite[Lemma p.~177]{Rud87}).

\begin{lemma}\label{lemtechnique}
Let $a$, $b$ be real numbers such that $a<b$, and $f:[a,b]\rightarrow\R$ 
be an absolutely continuous function. Let $J$ be a subset of $\{t\in [a,b]\ \vert\ f(t)=0\}$ os positive Lebesgue measure. Then $f'(t) = 0$ a.e. on $J$.
\end{lemma}

Using this lemma, and extremal equations \eqref{systemeadjoint}, one gets, for a.e. $t\in I$,
\begin{equation}\label{eqpsi}
\dot \Psi(t) = \frac{bC}{m(t)^2} \left( \|p_v(t)\| \|u(t)\| 
- \langle p_v(t),u(t)\rangle\right) + \Xi(r(t),v(t),m(t),p(t)) =0 ,
\end{equation}
where the function
\begin{equation*}
\begin{split}
\Xi(r,v,m,p) =  
 \frac{ Db}{m^2\|v\|} \langle p_v,v\rangle 
+
\frac{C}{m\|p_v\|} \bigg(  & \langle p_v,p_r\rangle 
+ \frac{\langle p_v, v\rangle}{m\Vert v\Vert}\left\langle \frac{\partial D}{\partial v},p_v\right\rangle 
\\ & 
+ \frac{\partial D}{\partial m} \frac{ \|p_v\|^2}{m\|v\|} 
- \frac{ D}{m}\frac{\langle p_v,v\rangle^2}{\|v\|^3} \bigg)
\end{split}
\end{equation*}
does not depend on $u$.
From Lemma \ref{lembang}, the relation \eqref{ucolinpv} holds almost everywhere, and hence the first term of \eqref{eqpsi} vanishes. Therefore,
\begin{equation}\label{dpsitdt2}
\dot \Psi(t) = \Xi(r(t),v(t),m(t),p(t))=0,
\end{equation}
for almost every $t\in I$ (actually over every subinterval of positive measure, since the above expression is continuous).

Relations \eqref{rel1} and \eqref{eqpsi} are two constraint equations, necessary for the existence of a singular arc.
Derivating once more, using Lemma \ref{lemtechnique}, leads to
\begin{equation}\label{psipp}
\ddot \Psi(t) = 0,\quad \textrm{a.e. on}\ I.
\end{equation}
The control $u$ is expected to appear explicitly in this latter relation. However, since calculations are too lengthy to be reported here, we next explain how \eqref{psipp} permits to derive an expression for $\alpha(t)$, and hence, from (\ref{ucolinpv}), an expression for $u(t)$.
When derivating (\ref{dpsitdt2}), the terms where the control $u$ appears are the terms containing $\dot{v}$, $\dot{p}_m$, and $\dot{m}$. Recall that $\dot{m}=-b\Vert u\Vert$, that $\dot{p}_m = \frac{1}{m} \langle p_v, -\frac{D(r,v)}{m}\frac{v}{\Vert v\Vert} + C\frac{u}{m} \rangle$, and that $\dot{v}$ is affine in $u$. Hence, since $\alpha(t)\geq 0$, it is not difficult to see that this derivation leads to an equation of the form
\begin{equation}\label{relfinale}
A(r,v,m,p_r,p_v,p_m)\alpha=B(r,v,m,p_r,p_v,p_m),
\end{equation}
almost everywhere on $I$.
This relation should be "generically" nontrivial, that is, the coefficient $A$ should not be equal to zero. This fact proves to hold true on numerical simulations. We explain below rigorously why this is true generically at least in the case of a scalar control (recall that we deal here with a three-dimensional control). For a scalar control, the control system (\ref{system}) is of the form
\begin{equation}\label{sysmono}
\dot{q}=f_0(q)+uf_1(q),
\end{equation}
where $f_0$ and $f_1$ are smooth vector fields, and $q$ is the state. In this case, it is well known (see e.g.\ \cite{BonChy03}) that, if $u$ is a singular control on $I$, then there must exist an adjoint vector $p$ such that
\begin{eqnarray}
\langle p,f_1(q)\rangle &=& 0\quad \textrm{on}\ I,  \label{r1} \\
\langle p,[f_0,f_1(q)]\rangle &=& 0\quad \textrm{on}\ I,   \label{r2}  \\
\langle p,[f_0,[f_0,f_1(q)]\rangle + u \langle p,[f_1,[f_0,f_1(q)]\rangle &=& 0\quad \textrm{a.e.\ on}\ I.  \label{r3}
\end{eqnarray}
The situation encountered here for 3D Goddard's problem is similar to that case: Equations (\ref{r1}), (\ref{r2}), (\ref{r3}), are respectively similar to Equations (\ref{rel1}), (\ref{dpsitdt2}), (\ref{psipp});
Equations (\ref{r1}), (\ref{r2}) (similarly, Equations (\ref{rel1}), (\ref{dpsitdt2})) are constraint equations, and Equation (\ref{r3}) (similarly, Equation (\ref{psipp})) permits in general to derive an expression for the control $u$. The vocable "generic" employed above can now be made more precise: it is proved in \cite{BonKup97} that there exists an open and dense (in the sens of Whitney) subset $\mathcal{G}$ of the set of couples of smooth vector fields such that, for every control system (\ref{sysmono}) with $(f_0,f_1)\in\mathcal{G}$,
the set where $\langle p,[f_1,[f_0,f_1(q)]\rangle$ vanishes has measure zero, and hence
Equation (\ref{r3}) always permits to derive $u$.
Additionaly, we can notice that the classical one-dimensional Goddard problem can be formulated as a particular case of the general 3D problem described here. In this case, it is well known that the second derivative of the switching function provides the expression of the singular control, so we can safely assume that \ref{relfinale} is nontrivial for the restriction to the 1D problem.
Based on these arguments, we should expect the coefficient $A$ of Equation (\ref{relfinale}) to be non zero in general. 
This is indeed the case in our numerical simulations presented next.
Of course, once $\alpha(t)$ has been determined, one has to check (numerically) that $0\leq \alpha(t)\leq 1$, so that the constraint $\Vert u\Vert\leq 1$ is indeed satisfied.
Here also, numerical simulations show the existence and admissibility of such singular arcs (see Section \ref{secnum}).


\subsection{Conclusion}
We sum up the previous results in the following theorem.

\begin{theorem}\label{thmcontroles}
Consider the optimal control problem of maximizing a final cost $g^0(t_f,x(t_f))$,
for the control system (\ref{system}), with initial conditions (\ref{initialconditions}) and 
final conditions $x(t_f)\in M_1$. We assume that Assumption $(H)$ holds. 
Let $u$ be an optimal control defined on $[0,t_f]$, associated to the trajectory $(r(\cdot),v(\cdot),m(\cdot))$. Then, there exist absolutely continuous mappings $p_r(\cdot):[0,t_f]\rightarrow\R^3$, $p_v(\cdot):[0,t_f]\rightarrow\R^3$, $p_m(\cdot):[0,t_f]\rightarrow\R$, and a real number $p^0\leq 0$, such that $(p_x(\cdot),p_v(\cdot),p_m(\cdot),p^0)$ is nontrivial, and such that Equations (\ref{systemeadjoint}) hold a.e.\ on $[0,t_f]$.
Define the switching function $\Psi$ on $[0,t_f]$ by
$$\Psi(t) = \frac{C}{m(t)}\Vert  p_v(t)\Vert -bp_m(t) .$$
Then,
\begin{itemize}
\item
if $\Psi(t) < 0$ then $u(t)=0$;
\item
if $\Psi(t) > 0$ then $u(t)=\frac{p_v(t)}{\Vert p_v(t)\Vert}$;
\item
if $\Psi(t)=0$ on a subset $I\subset[0,t_f]$ of positive Lebesgue measure, then Equation (\ref{dpsitdt2}) must hold on $I$, and
$$u(t)=\alpha(t)\frac{p_v(t)}{\Vert p_v(t)\Vert}\quad\textrm{a.e.\ on}\ I,$$
where $\alpha(t)\in[0,1]$ is determined by (\ref{relfinale}).
\end{itemize}
\end{theorem}

\begin{remark}
The optimal control is piecewise either equal to zero, or saturating the constraint with the direction of $p_v(t)$, or is singular. Notice that, in all cases, it is collinear to $p_v(t)$, with the same direction.
\end{remark}

\begin{remark}[Optimality status]
The maximum principle is a necessary condition for optimality. Second-order sufficient conditions are usually characterized in terms of conjugate points (see e.g. \cite{AgSachkov,BCT}. Unfortunately standard theories do not apply here for two reasons: first, the equation in $m(t)$ involves the term $\Vert u(t)\Vert$ which is not smooth; second, the structure of trajectories stated in the theorem involves both bang arcs and singular arcs, and up to now a theory of conjugate points that would treat this kind of trajectory.

We mention however below a trick, specific to the form of our system, which permits to apply the standard theory of conjugate points on every subinterval $J$ of $[0,t_f]$ on which $u$ is singular and $0<\Vert u(t)\Vert<1$. Let $J$ be such a subinterval. Then, $\dot{m}\neq 0$ a.e.\ on $J$, and the system can be reparametrized by $-m(t)$. Then, denoting $q=(r,v)$, system \eqref{system} yields
$$\frac{dq}{dm} = \frac{1}{\Vert u\Vert}f(m,q)+\frac{u_1}{\Vert u\Vert}g_1(m,q)+\frac{u_2}{\Vert u\Vert}g_2(m,q)+\frac{u_3}{\Vert u\Vert}g_3(m,q).$$
Now, set
$$v=\frac{1}{\Vert u\Vert},\ \textrm{and}\ \frac{u_1}{\Vert u\Vert}=\cos\theta_1\cos\theta_2,\ \frac{u_2}{\Vert u\Vert}=\cos\theta_1\sin\theta_2,\ \frac{u_3}{\Vert u\Vert}=\sin\theta_2,$$
and consider as new control the control $\tilde u=(v,\theta_1,\theta_2)$. Notice that the controls $\theta_1$ and $\theta_2$ are unconstrained, and that $v$ must satisfy the constraint $v\geq 1$. However, along the interval $J$ it is assumed that $0<\Vert u(t)\Vert<1$, and thus $v$ does not saturate the constraint. Hence, the standard theory of conjugate points applies and the local optimality status of the trajectory between its extremities on $J$ can be numerically checked, for instance using the code \textit{COTCOT} (Conditions of Order Two and COnjugate times), available on the web\footnote{http://www.n7.fr/apo/cotcot},
developed in \cite{BCT}. This reference provides algorithms to compute the first
conjugate time (where the trajectory ceases to be optimal) along a smooth extremal curve, based on theoretical developments of geometric optimal control using second order optimality conditions. The computations are related to a test of positivity of the intrinsic second order derivative or a test of singularity of the extremal flow.

It can be checked as well that every smooth sub-arc of the trajectory is locally optimal between its extremities. However, the problem of proving that the \textit{whole} trajectory (i.e., a succession of bang and singular arcs) is locally optimal is open. Up to now no conjugate point theory exists to handle that type of problem. Of course, one could make vary the times of switchings but this only permits to compare the trajectory with other trajectories having exactly the same structure. A sensitivity analysis is actually required to treat trajectories involving singular subarcs.
\end{remark}


\section{Numerical experiments}\label{secnum}

In this section, we provide numerical simulations showing the relevance of singular arcs in the complete Goddard's Problem.
For given boundary conditions, the optimal trajectory is first computed using \textit{indirect methods} (shooting algorithm) combined with an homotopic approach.
Then we use a \textit{direct method} (based on the discretization of the problem) to check the obtained solution. 
All numerical experiments were led on a standard computer (Pentium 4, 2.6 GHz).

\subsection{Numerical values of the parameters of the model}
We implement the optimal control problem of maximizing $m(t_f)$ for the system \eqref{system}, with the constraint \eqref{constraint}.
The equations of motion can be normalized with respect to $r(0)$, $m(0)$, and $g_0$. We follow \cite{Obe90} (in which 2D-trajectories with maximization of the final velocity are studied), and set the following parameters. 
\begin{itemize}
\item The distance unit is the Earth radius $R_T=6378\ 10^3$ m.
\item Maximal thrust modulus $C=3.5$; $b=7$.
\item Gravity $g(r) = \frac{g_0}{\|r\|^3} r $, with $g_0 = 1$.
\item Drag $D(r,v) = K_D \|v\|^2 e^{-500(\|r\|-1)}$ with $K_D=310$.
\item Initial and final conditions
$$
\begin{array}{l}
r_0=(0.999949994\ 0.0001\ 0.01), \quad v_0=(0\ 0\ 0), \quad m_0=1,\\
r_f=(1.01\ 0\ 0), \quad v_f\  \textrm{is free},\quad  m_f\  \textrm{is free.}\\
t_f\  \textrm{is free.}
\end{array}
$$
\end{itemize}


\subsection{Numerical simulations with indirect methods}
In our simulations presented hereafter, we prefer to express the objective of the optimal control problem in the following form.\\
Maximizing $m(t_f)$ is equivalent to minimizing the cost
$$\int_0^{t_f} \Vert u(t)\Vert dt,$$
and we assume that there are no minimizing abnormal extremals, therefore the adjoint vector can be normalized so that $p^0=-1$. 
The results of the simulations are consistent with this assumption.\\

According to Section \ref{secPMP}, the Hamiltonian of the optimal control problem under consideration is
$$
H=\langle p_r,v\rangle + \left\langle p_v, -\frac{D(r,v)}{m}\frac{v}{\Vert v\Vert} - g(r) + C\frac{u}{m}  \right\rangle -(1 + b p_m) \Vert u\Vert,
$$
The only difference with the Hamiltonian in \ref{secPMP} for the $Max\ m(t_f)$ objective is the additional ``$-1$'' in the $\|u\|$ term, which leads to the switching function $\psi(t) = \frac{C}{m(t)}\Vert  p_v(t)\Vert - (1 + b p_m(t)) $,
\begin{itemize}
\item
if $\psi(t) < 0$ then $u(t)=0$;
\item
if $\psi(t) > 0$ then $u(t)=\frac{p_v(t)}{\Vert p_v(t)\Vert}$;
\item
if $\psi(t)=0$ on $I\subset[0,t_f]$, then Equation (\ref{dpsitdt2}) must hold on $I$, the control $u$ is singular, and
$$u(t)=\alpha(t)\frac{p_v(t)}{\Vert p_v(t)\Vert}\quad\textrm{a.e.\ on}\ I,$$
where $\alpha(t)\in[0,1]$ is determined by (\ref{relfinale}). We check numerically that $0 \leq \alpha(t)\leq 1$.
\end{itemize}

Furthermore, on a singular subarc, derivating the switching function twice yields the expression of $\alpha$ via a relation of the form $A(x,p)\alpha=B(x,p)$, see \ref{relfinale}.
The computations are actually quite tedious to do by hand, and we used the symbolic calculus tool \textsc{Maple}.
The expressions of $A$ and $B$ are quite complicated and are not reported here.\\

The free final time problem is formulated as a fixed final time one via the usual time transformation $t = t_f\ s$, with $s \in [0,1]$ and $t_f$ an additional component of the state vector, such that $\dot{t_f}=0$ and $t_f(0), t_f(1)$ are free, with the associated costate satisfying $\dot{p_{t_f}} = -H$.\\

Transversality conditions on the adjoint vector yield $p_v(1)=(0\ 0\ 0)$, $p_m(1)=0$, and $p_{t_f}(0) = p_{t_f}(1) = 0$.\\

\subsubsection{Homotopic approach}
In the indirect approach, it is necessary to get some information on the structure of the solutions, namely, to know a priori the number and approximate location of singular arcs. To this aim, we perform a continuation (or \emph{homotopic}) approach, and regularize the original problem by adding a quadratic ($\|u\|^2$) term to the objective, as done for instance in \cite{Martinon05,Schwartz96}.
The general meaning of continuation is to solve a difficult problem by starting from the known solution of a somewhat related, but easier problem.
By related we mean here that there must exist a certain application $h$, called a \emph{homotopy}, connecting both problems.
Here, we regularize the cost function by considering an homotopic connection with an energy, 
\begin{equation}\label{coutcontin}
\int_0^{t_f} \ (\|u(t)\| + (1 - \lambda) \|u(t)\|^2) \ dt,
\end{equation}
where the parameter of the homotopy is $\lambda \in [0,1]$.
The resulting perturbed problem $(P_\lambda)$ has a strongly convex Hamiltonian (with respect to $u$), with a continuous optimal control, and is much easier to solve than $(P)=(P_1)$. 
Assuming we have found a solution of $(P_0)$, we want to follow the zero path of the homotopy $h$ until $\lambda=1$, in order to obtain an approximate solution of $(P)$ (or at least sufficient information).
The continuation can be conducted manually, by finding a suitable sequence $(\lambda_k)$ from $0$ to $1$.
However, finding such a sequence can be quite difficult in practice, which is why we chose here to perform a full path-following continuation.  
Extensive documentation about path following methods can be found in \cite{AlGe90}. We use here a piecewise-linear (or \emph{simplicial}) method, whose principle is recalled briefly below. 
The reason behind the choice of this method over a more classical predictor-corrector continuation (such as detailed for instance in \cite{Deu04}) is that we expect the problem to be ill-conditioned, due to the presence of singular arcs, which is indeed the case in the numerical experiments.

\medskip

\noindent{\bf Simplicial methods.}
PL continuation methods actually follow the zero path of the homotopy $h:\mathbf{R}^{n+1} \to \mathbf{R}^n$ by building a piecewise linear approximation of $h$.
The search space $\mathbf{R}^{n+1}$ is subdivided into cells, most often in a particular way called \emph{triangulation} in \emph{simplices}.
This is why PL continuation methods are often referred to as simplicial methods.
The main advantage of this approach is that it imposes extremely low requirements on the homotopy $h$: since no derivatives are used, continuity is in particular sufficient, and should not even be necessary in all cases.

\begin{definition}[Simplices and faces] A simplex is the convex hull of $n+1$ affinely independent points (called the vertices) in $\mathbf{R}^n$, while a $k$-face of a simplex is the convex hull of $k$ vertices of the simplex ($k$ is typically omitted for $n$-faces, which are just called faces, or facets).
\end{definition}

\begin{definition}[Triangulation] A triangulation is a countable family $T$ of simplices of $\mathbf{R}^n$ such that
the intersection of two simplices of $T$ is either a face or empty, and such that
$T$ is locally finite (a compact subset of $\mathbf{R}^{n}$ meets finitely many simplices).
\end{definition}

\begin{figure}[h]
\begin{center}
\scalebox{0.3}[0.3]{\includegraphics*{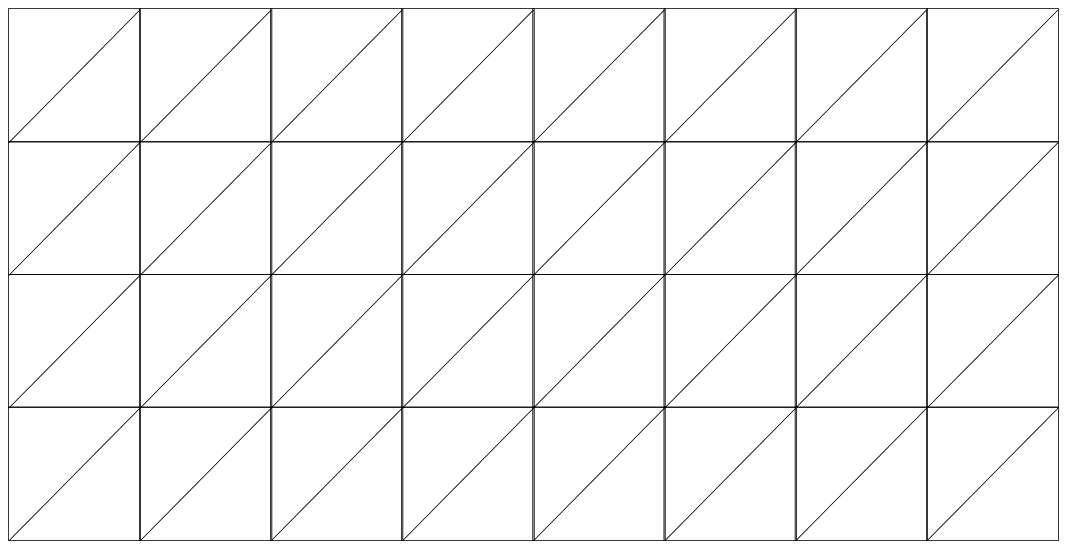}}
\hspace*{1cm}
\scalebox{0.32}[0.32]{\includegraphics*{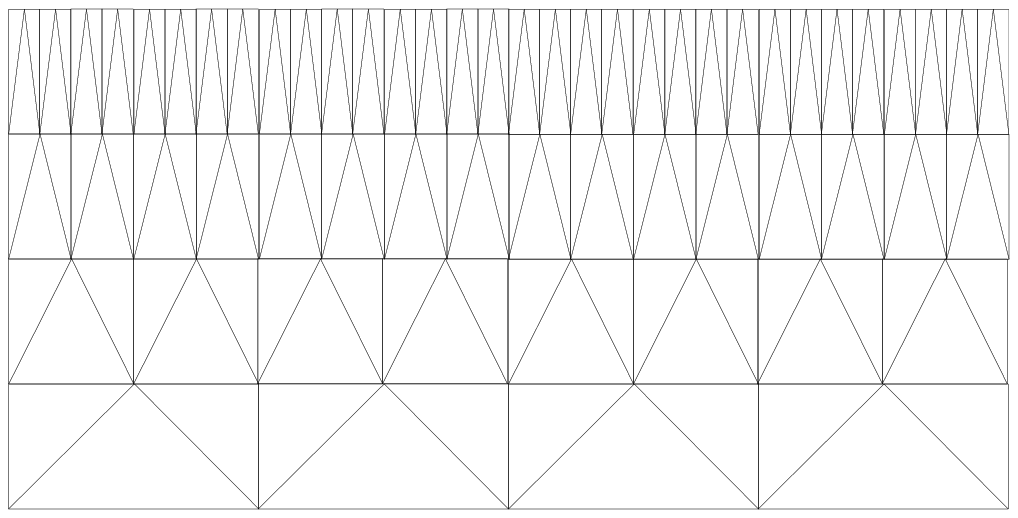}}
\hspace*{1cm}
\scalebox{0.32}[0.32]{\includegraphics*{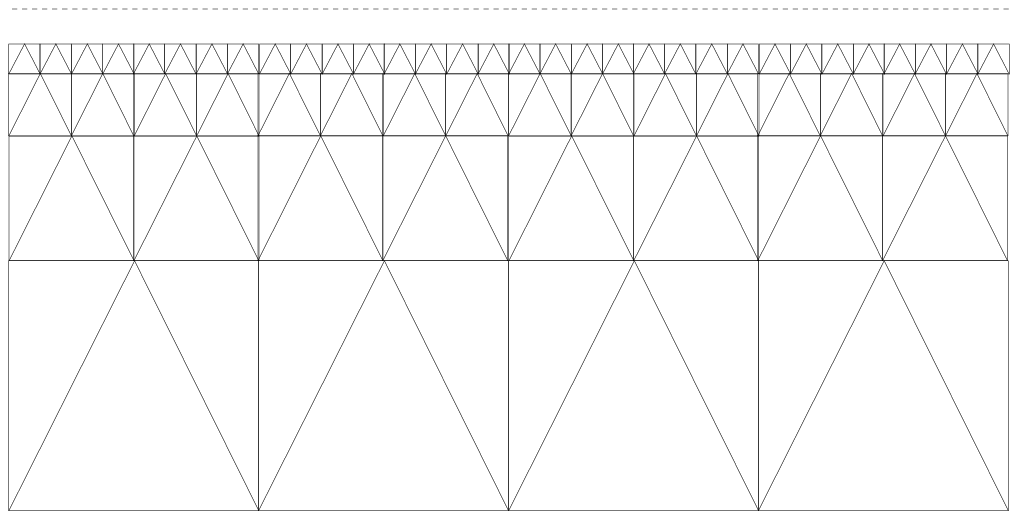}}
\caption{Illustration of some well known triangulations of $\mathbf{R} \times [0,1]$ ($[0,1[$ for $J_3$):
Freudenthal's uniform $K_1$,  and Todd's refining $J_4$ and $J_3$}
\end{center}
\end{figure}

\begin{definition}[Labeling] We call labeling a map $l$ that associates a value to the vertices $v_i$ of a simplex.
We label here the simplices by the homotopy $h$: $l(v^i)=h(z^i,\lambda^i)$, where $v^i=(z^i,\lambda^i)$.
Affine interpolation on the vertices thus gives a PL approximation $h_T$ of $h$.
\end{definition}

\begin{definition}[Completely labeled face] A face $[v_1,..,v_n]$ of a simplex is said completely labeled if and only if it contains a solution $v_\epsilon$ of the equation $h_T(v) = \vec{\epsilon}=(\epsilon,..,\epsilon^n)$, for every $\epsilon >0$ sufficiently small.
\end{definition}

\begin{lemma}[{\cite[Chapter 12.4]{AlGe90}}] 
Each simplex possesses either zero or exactly two completely labeled faces (called a \emph{transverse} simplex in the latter case).
\end{lemma}

The constructive  proof of this property, which gives the other completely labeled face of a simplex that already has a known one, is often referred to as \emph{PL step, linear programming step}, or \emph{lexicographic minimization}.
Then there exists a unique transverse simplex sharing this second completely labeled face, that can be determined via the \emph{pivoting rules} of the triangulation.\\

A simplicial algorithm thus basically follows a sequence of transverse simplices, from a given first transverse simplex with a completely labeled face at $\lambda=0$, to a final simplex with a completely labeled face at $\lambda=1$ (or $1-\epsilon$ for refining triangulations that never reach 1), which contains an approximate solution of $h(z,1)=0$.

\subsubsection{Preliminary continuation on the atmosphere density}
In our case, even solving the regularized problem $(P_0)$ is not obvious, due to the aerodynamic forces (drag).
For this reason, we introduce a preliminary continuation on the atmosphere density, starting from a problem without atmosphere.
Technically, this is done by using an homotopy with the modified parameter
$$
K_D^\theta = \theta K_D, \quad \theta \in [0,1],
$$
where $K_D$ appears in the model of the drag.
The shooting method for the problem without atmosphere at $\theta=0$ converges immediately with the trivial starting point $z_0 =(0.1\ 0.1\ 0.1\ 0.1\ 0.1\ 0.1\ 0.1)$.
We would like to emphasize the fact that we have here no difficulties to find a starting point for the shooting method.
The path following is then achieved with an extremely rough integration formula (Euler with only 25 steps), since we just seek a starting point for the main homotopy.
Thanks to the robustness of the simplicial method, we can afford such a low precision to save computational time. 
The border at $\theta=1$ is reached after crossing about $120\ 000$ simplices, for a CPU time of 48 seconds.\\

\begin{figure}[h]
\begin{center}
\includegraphics[width=12cm,height=6cm]{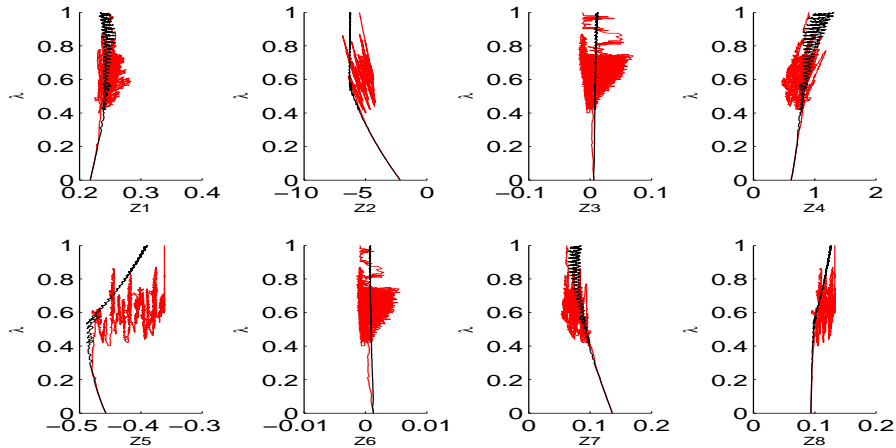}
\end{center}
\caption{Path following for the atmosphere homotopy: fixed uniform triangulation (grey) and adaptive meshsize (black).}\label{figpath}
\end{figure}

\begin{remark}
The adaptive meshsize algorithm described in \cite{Martinon05} here strongly reduces the oscillations along the zero path, as shown on Figure \ref{figpath}, which decreases the number of simplices required to reach $\theta=1$.
We can see that the path following using a fixed uniform meshsize actually converges to another point, which corresponds to an incorrect solution (the final condition on $r_2$ is not satisfied).
\end{remark}

\begin{figure}[h]
\begin{center}
\includegraphics[width=8cm,height=7cm]{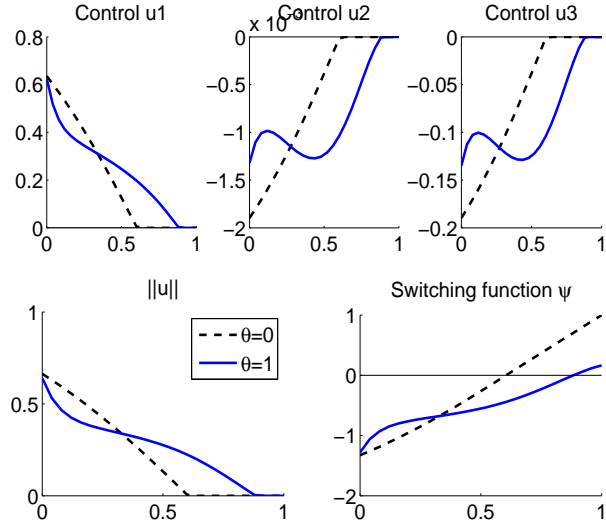}
\end{center}
\caption{Regularized problem $(P_0)$: solutions for $\theta=0$ (no atmosphere) and $\theta=1$ (normal atmosphere).}\label{figatm}
\end{figure}

The solution we obtain is sufficient to initialize the shooting method at the beginning of the main homotopy.
Figure \ref{figatm} represents the solutions of the regularized problem $(P_0)$ for $\theta=0$ and $\theta=1$, i.e., without atmosphere and with a normal atmosphere.\\

Notice that a direct continuation on the atmosphere with the original non regularized problem $(P)$ fails. During the continuation, the process abruptly diverges at a certain value for $\theta$, certainly due to the appearance of the singular arc.

\subsubsection{Main continuation on the quadratic regularization}
\begin{figure}[h]
\begin{center}
\includegraphics[width=10cm,height=8cm]{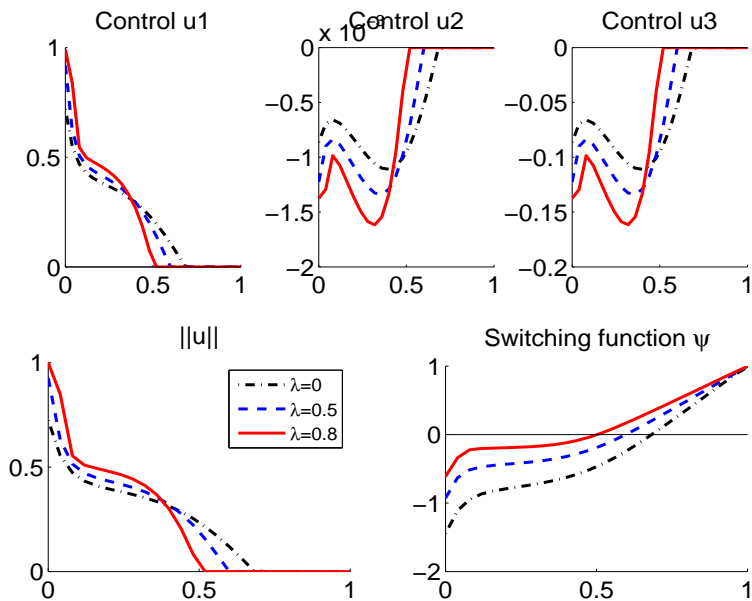}
\caption{Main homotopy - solutions for $\lambda=0,0.5$ and $0.8$.}\label{figcontin}
\end{center}
\end{figure}

We now perform the main continuation on the cost \eqref{coutcontin}.
Figure \ref{figcontin} represents the solutions for $\lambda=0,0.5$ and $0.8$. It is visible that this continuation process permits to detect the singular structure of the solution.
The shape on the switching function and of the control norm graphs are particularly interesting concerning suspicion of singular arcs. Indeed,
we observe that, on a certain time interval (roughly $[0.1, 0.4]$), the switching function comes closer to zero as $\lambda$ increases, while the control norm keeps values in $(0,1)$. Along the solution for $\lambda=0.8$, we can guess the appearance of a small arc where $\|u\|=1$ at the beginning.
These facts strongly suggest the appearance of a singular arc.

With a fixed meshsize of $10^{-4}$, the path following takes about 900000 simplices and 350 seconds to reach $\lambda=0.8$, again with an extremely rough integration (Euler, 25 steps).
Trying to go further becomes extremely difficult since we lose the singular structure and encounter trajectories with incorrect bang-bang structures. However the knowledge of the solution for $\lambda=0.8$ happens to be sufficient to solve the problem: it provides a good starting point for which the shooting method applied to the original problem $(P)$ converges.\\

\begin{remark} 
This path is more difficult to follow than the previous one for the atmosphere homotopy, and the adaptive meshsize algorithm does not work well. We thus use a fixed meshsize to perform this homotopy.
\end{remark}

\subsubsection{Shooting method applied to the original problem $(P)$}
When implementing a shooting method (see for instance \cite{Bet01,Bul71,Hilt90,StoBul93}),
the structure of the trajectory has to be known a priori. The structure of the control must be prescribed here by assigning a fixed number of interior switching times that correspond to junctions between nonsingular and singular arcs.
These times $(t_i)_{i=1..n_{switch}}$ are part of the shooting unknowns and must satisfy some switching conditions.
Each arc is integrated separately, and matching conditions must be verified at the switching times, as drawn on the diagram below.

\begin{quote}

Unknown: $z$

\centerline{\footnotesize
\begin{tabular}{|c|c|c|c|c|c|c|}
\hline
IVP unknown at $t_0$ & $(x^1,p^1)$ & ... & $(x^s,p^s)$ & $t_1$ & ... & $t_s$\\
\hline
\end{tabular}
}

Value: $S_{Sing}(z)$

\centerline
{\footnotesize
\begin{tabular}{|c|c|c|c|c|c|}
\hline
$Switch_{cond}(t_1)$ & $Match_{cond}(t_1)$ & ...  & $Switch_{cond}(t_s)$ & $Match_{cond}(t_s)$ & TC($t_f$)  \\
\hline
\end{tabular}
}
\end{quote}

Here, matching conditions reduce to imposing state and costate continuity at the switching times.

A switching condition indicates a change of structure, which corresponds here to an extremity of a singular arc. Along such a singular arc, it is required that $\psi=\dot\psi=0$. The control is computed using the relation $\ddot{\psi}=0$. Therefore, using this expression of the control, switching conditions consist in imposing either $\psi=0$ at the extremities of the singular arc, or $\psi=\dot{\psi}=0$ at the beginning of the arc. In our simulations, we choose the latter solution which happens to provide better and more stable results.

\medskip

The previous results, obtained with an homotopic approach, provide an indication on the expected structure of the optimal trajectory for the original problem $(P)$. Inspection of Figure \ref{figcontin} suggests to seek a solution involving a singular arc on an interval $[t_1,t_2]$, with $t_0<t_1<t_2<t_f$. As a starting point of the shooting method, we use the solution previously obtained with the homotopy on the cost at $\lambda=0.8$.\\

The IVP integration is performed with the \textsc{radau5} code (see \cite{Hairer03}), with absolute and relative tolerances of, respectively, $10^{-6}$ and $10^{-6}$. The shooting method converges in 17 seconds, with a shooting function of norm $5\ 10^{-4}$.
The condition number for the shooting function is quite high (about $10^{12}$), which was expected. 
The overall execution time of the whole approach (preliminary atmosphere homotopy, regularization homotopy, final shooting) is about 400 seconds.\\

\begin{figure}[h]
\centerline{\includegraphics[width=4cm,height=4cm]{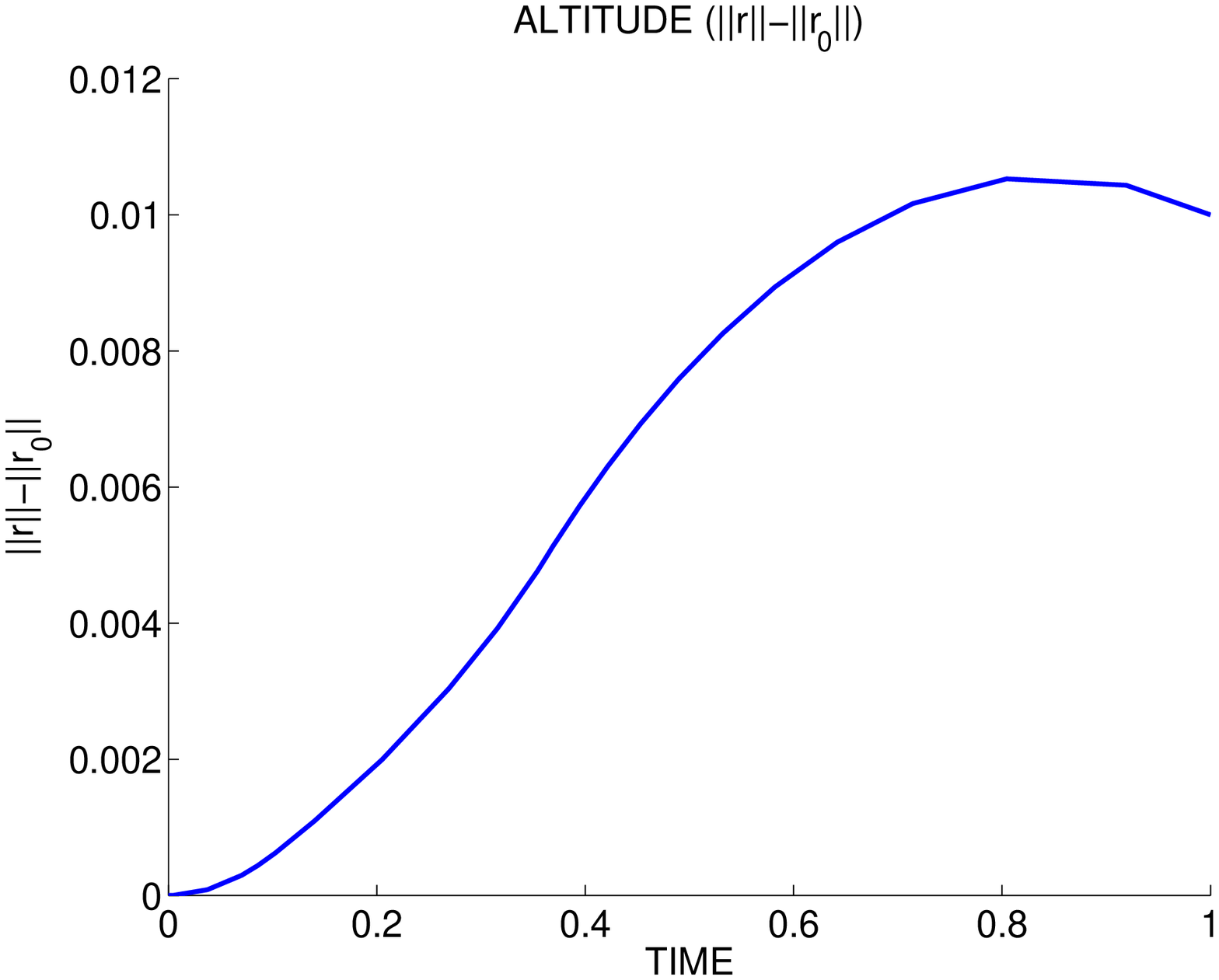} 
\includegraphics[width=4cm,height=4cm]{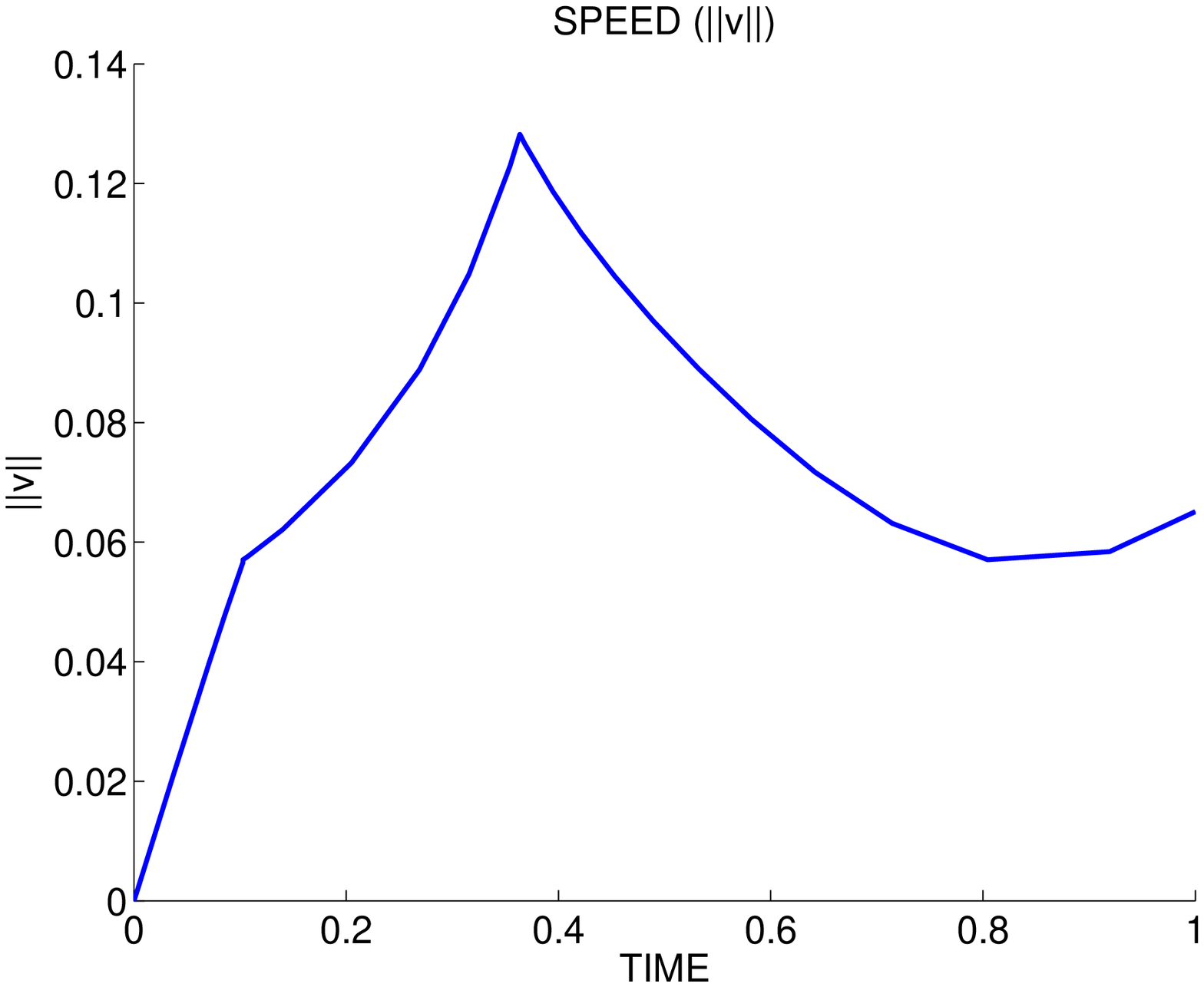} 
\includegraphics[width=4cm,height=4cm]{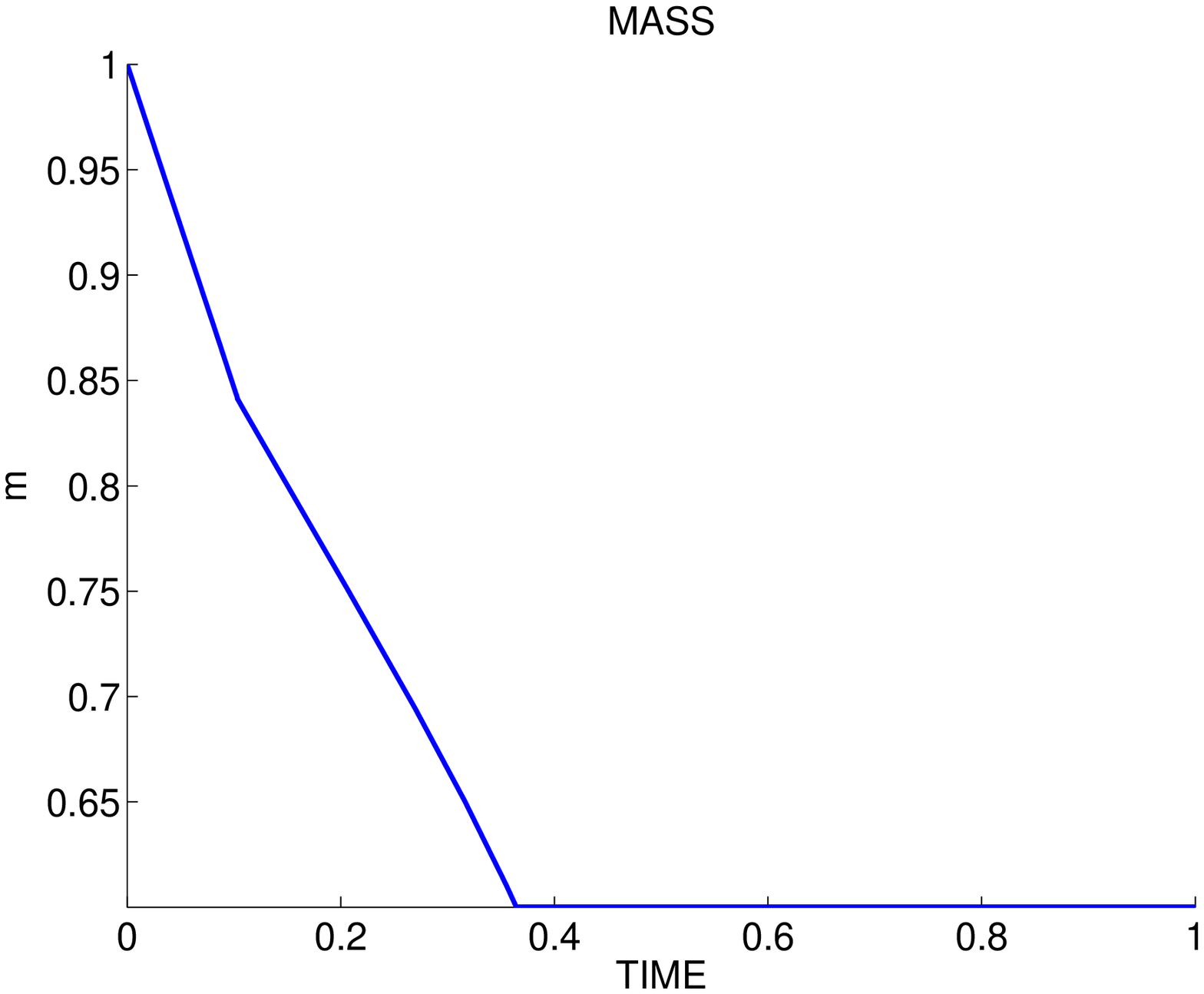}}
\caption{Solution with singular arc: altitude, speed and mass.} \label{figsol}
\end{figure}

At the solution, the free final is $0.2189$, and the objective value is $0.3994$, which corresponds to a final mass of $0.6006$.
The evolution of altitude, speed and mass during the flight are represented on Figure \ref{figsol}.

\begin{figure}[h]
\centerline{\includegraphics[width=10cm,height=8cm]{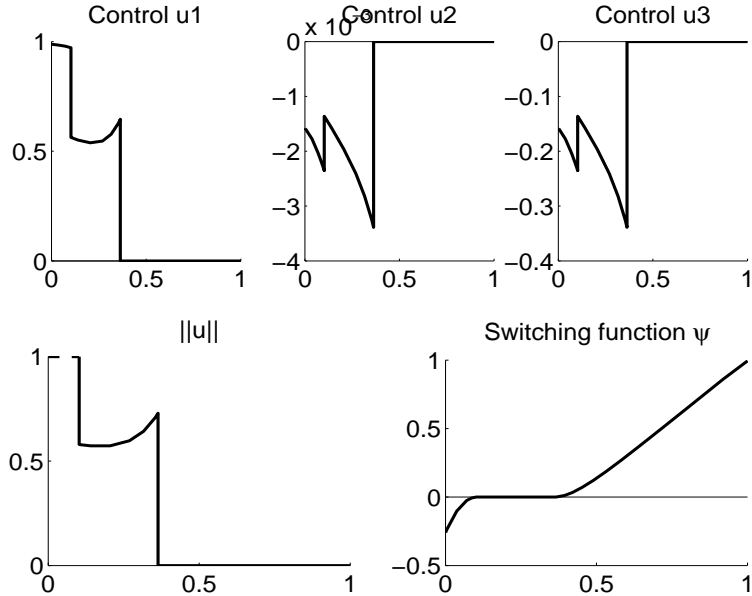}}
\caption{Solution with singular arc: control and switching function.} \label{figcontrol}
\end{figure}

We show on Figure \ref{figcontrol} the control and switching function.
The singular arc is clearly visible on the control norm graph.

\subsection{Numerical simulations with direct methods}
In order to validate the solution obtained previously with the shooting algorithm, we next implement a direct method.
Although direct methods can be very sophisticated (see for instance \cite{Bet01,WachBieg06}), we here use a very rough formulation, since our aim is just to check if the results are consistent with our solution.
We discretize the control using piecewise constant functions, and the state is integrated on $[0,t_f]$ with a basic fixed step Runge-Kutta fourth order formula.
The values of the control at the discretization nodes, as well as the final time $t_f$, thus become the unknowns of a nonlinear constrained optimization problem, the constraints being the final conditions for the state.
To solve the optimization problem, we use the \textsc{ipopt} solver, which implements an interior point algorithm with a filter line-search method (see \cite{WachBieg06} for a complete description).

\begin{figure}[h]
\centerline{\includegraphics[width=10cm,height=8cm]{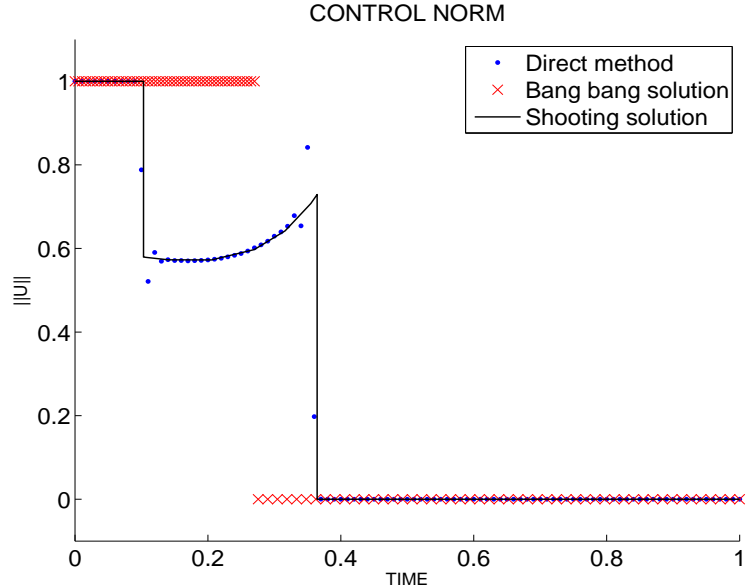}}
\caption{Control norm for the shooting and direct method.} \label{figdirecte}
\end{figure}

With standard options, the algorithm converges after 193 iterations (and 210 seconds) to a solution with a final time of $0.2189$ and a criterion value of $0.3997$.
This solution is clearly consistent with the results of the shooting method, as shown on Figure \ref{figdirecte}, which represents the norm of the control for the shooting method solution, the direct method solution, and a bang-bang reference solution (see below).

\subsubsection*{Comparison with a bang-bang solution}
Recall that the usual launch strategy consists in implementing piecewise controls either saturating the constraint or equal to zero. To prove the relevance of the use of singular controls in the control strategy, we next modify slightly the formulation above in order to find a bang-bang solution. Our aim is to demonstrate that taking into account singular arcs in 
the control strategy actually improves (as expected) the optimization criterion.

We implement a \emph{``on-off''} structure, with only one switching time $t_\textrm{off}$. The control is chosen so as to satisfy $\Vert u(t)\Vert=1$ for $t_0<t<t_\textrm{off}$, and $u(t)=0$ for $t_\textrm{off}<t<t_f$.
Here, the unknowns of the optimization problem are $t_f$, $t_\textrm{off}$ and the direction of the control at the discretization nodes before $t_\textrm{off}$.
We obtain a solution with $t_f=0.2105$, $t_\textrm{off}=0.0580$, and the value of the criterion is $0.4061$, which represents a loss of about $1.6\%$ compared to the solution with a singular arc.
On this academic example, the gain of the optimal strategy, involving a singular arc, over a pure bang-bang strategy, is quite small.
This simplified problem is a first step in the study of a realistic launcher problem, and permits to illustrate the method.

\bibliography{hjb}
\bibliographystyle{plain}

\tableofcontents

\end{document}